\documentclass[11pt,a4paper]{article}
\usepackage{amsmath}
\usepackage{amsfonts}
 \usepackage{amsthm}
 \usepackage{hyperref}
 \usepackage[usenames]{color} 
 \usepackage{mathrsfs}
\theoremstyle{plain}
\newtheorem{thm}{Theorem}[section]

\newtheorem{prop}[thm]{Proposition}
\newtheorem{wn}[thm]{Corollary}
\theoremstyle{definition}
\newtheorem{rem}[thm]{Remark}
\newcommand{\proofof}[1]{\vskip 0pt\noindent\textbf{\textit{Proof of #1. }}}

\newcommand{\R}{\mathbb R}
\newcommand{\N}{\mathbb N}

\newcommand{\C}{\mathbb C}

\newcommand{\abs}[1]{\left\vert#1\right\vert} 

\newcommand{\ind}[1]{1\mkern-7mu1_{\{#1\}}}
\newcommand{\indn}{1\mkern-7mu1}

\numberwithin{equation}{section}
 \title{\normalsize{\textbf{Bifractional Brownian motion for $H>1$ and $2HK\le  1$ %and decomposition of fractional Brownian motion with small Hurst parameter.
}}}
 
  \author{
 \small{ANNA TALARCZYK}\thanks{Institute of Mathematics, University of Warsaw, ul. Banacha 2, 02-097 Warsaw, Poland,  \hbox{e-mail:}
annatal@mimuw.edu.pl. Research supported in part by  National Science Centre (Narodowe Centrum Nauki), Poland, grant  2016/23/B/ST1/00492.} \\
}
\date{}

\begin{document}
 \maketitle
  \begin{abstract} 
   Bifractional Brownian motion on $\mathbb{R}_+$ is a two parameter centered Gaussian process with covariance function:
$$
   R_{H,K} (t,s)=\frac 1{2^K}\left(\left(t^{2H}+s^{2H}\right)^K-\abs{t-s}^{2HK}\right), \qquad s,t\ge 0.
$$
This process has been originally introduced by Houdr\'e and Villa in \cite{HoudreVilla} for the range of parameters $H\in (0,1]$ and $K\in (0,1]$. Since then, the range of parameters, for which $R_{H,K}$ is known to be nonnegative definite has been somewhat extended, but the full range is still not known.
  We give an elementary proof that $R_{H,K}$ is nonnegative definite for parameters $H,K$ satisfying $H>1$  and $0<2HK\le 1$.
  We  show that $R_{H,K}$ can be decomposed into a sum of two nonnegative definite functions.  As a side product we  obtain a decomposition of the fractional Brownian motion with Hurst parameter $H<\frac 12$ into a sum of time rescaled Brownian motion and  another independent self-similar Gaussian process.
  We also discuss some simple properties of bifractional Brownian motion with  $H>1$.

 \end{abstract}
 % \newpage
 
 \bigskip
 
 \section{Introduction and results}
 Bifractional Brownian motion (bifBm) has been first introduced by Houdr\'e and Villa in \cite{HoudreVilla} for parameters $0<H\le 1$, $0<K\le 1$ as a centered Gaussian process  with covariance function
 \begin{equation}
    R_{H,K} (t,s)=\frac 1{2^K}\left(\left(t^{2H}+{s}^{2H}\right)^K-\abs{t-s}^{2HK}\right), \qquad s,t\ge 0.
 \label{e:bifBm}
 \end{equation}
 In fact, for the same values of $H$ and $K$ in \cite{HoudreVilla} bifBm on the whole real line was studied, that is a centered Gaussian process on $\R$ with covariance function
 \begin{equation}
\frac 1{2^K}\left(\left(\abs t^{2H}+\abs{s}^{2H}\right)^K-\abs{t-s}^{2HK}\right), \qquad s,t\in\R.
\label{bifBmR}
\end{equation}
As noticed in \cite{Gorostiza}, for the same range of parameters the function \eqref{bifBmR}  appeared earlier in \cite{Berg}, Exercise 2.12 in Chapter 3.2 as an example of a positive definite kernel.

 In this article we will only discuss processes on the half-line, as our methods cannot be extended to the whole real line.
 
 If  $\xi$ is a centered Gaussian process with covariance function of the form \eqref{e:bifBm}, it is clearly self similar with index $HK$, that is, for any $a>0$ the process $(\xi_{at})_{t\ge 0}$ is equivalent in law to $(a^{HK}\xi_t)_{t\ge 0}$
 
 For $K=1$ bifractional Brownian motion reduces to the well known fractional Brownian motion (\cite{MandelbrotVanNess}). 
 Recall that fractional Brownian motion is a centered Gaussian process with covariance 
\begin{equation}
S_{H}(s,t)=\frac 12\left(s^{2H}+t^{2H}-\abs{t-s}^{2H}\right), \qquad s,t\ge 0.
\label{e:fBm}
\end{equation}
Up to a multiplicative constant, this is the only centered Gaussian process which is self similar and has stationary increments. 

Bifractional Brownian motion for $K\neq 1$ does not have stationary increments but it retains some of the properties of fractional Brownian motion with parameter $HK$, e.g. nonsemimartingale property, existence of modification with H\"older continuous trajectories with index $<HK$, nontrivial $1/HK$ variation.

Bifractional Brownian motion attracted quite a lot of attention. Various  properties of bifBm were studied, as well as stochastic integrals with respect to bifBm. By no means complete list includes for example: decompositions of bifBm and fractional Brownian motion \cite{LeiNualart}, \cite{Ma}, large time asymptotics \cite{MaejimaTudor}, sample path properties \cite{TudorXiao}, nonsemimartingale property also for $HK=\frac 12$ and introduction of stochastic integral \cite{RussoTudor}, relation to the stochastic heat equation and $p$-variation of bifractional Brownian motion \cite{Swanson}. Bifractional Brownian motion also often appears as an illustrative example for various more general theorems on self-similar processes e.g. \cite{Nualart} and others.

Bifractional Brownian motion appeared (admittedly, only with parameter $H=\frac 12$) in several natural models. For example  in
\cite{Swanson} in relation to stochastic heat equation, in \cite{Marouby} in a model of micropulses, or \cite{BGTnegative} in relation to branching particle systems.  In the latter paper the bifractional Brownian motion was obtained indirectly, as negative subfractional Brownian motion which in turn can be expressed as the integral of bifractional Brownian motion, see also \cite{extfBm}.

 \medskip
 
 Bardina and Es-Sebaiy \cite{BardinaEsSebaiy} showed that the range of parameters, for which one can define bifractional Brownian motion, can be extended to $0<H\le 1 $ and $0<K\le \min (2, \frac 1H)$. More recently, Durieu and Wang in \cite{DurieuWang} (see Remark 2.5 therein) showed that bifBm also exists if $H>1$ and $K=\frac 1{2H}$.
 
The exact range of parameters for which \eqref{e:bifBm} is a covariance function does not seem to be known. In  \cite{LifshitsVolkova} Lifshits and Volkova
study existence of bifBm in the case of $H>1$. 
A simple calculation shows that $K\le \frac 1H$ is 
a necessary condition for existence of bifBm. Otherwise
the inequality 
\begin{equation*}
 \abs{R_{H,K}(1,t)}^2\le R_{{H,K}}(1,1)R_{H,K} (t,t) 
\end{equation*}
would be violated for large $t$.
The authors proceed to a more subtle investigation, using the Lamperti transformation and  the property, that by self-similarity, the existence of a self similar process $\xi$ with covariance function \eqref{e:bifBm} is equivalent to existence of a stationary process 
$\zeta$, where $\zeta$ and $\xi$ are related by
\begin{equation*}
 \zeta(t)=e^{-HKt}\xi(e^t), \qquad t\in \R.
\end{equation*}
In \cite{LifshitsVolkova}  a  necessary  and sufficient condition is obtained in a non explicit form
(see Proposition 4.1 and formula (5.1) therein).  This condition is difficult to verify, since it relies on  a comparison of two quite complicated functions depending on parameters $H$ and $K$.
The authors only provide some numerical simulations showing for which values of $H$ and $K$ it is satisfied. 
A necessary condition in a different form is given by 
 Proposition 3.2 of the same paper. Again, it is not explicit.

\medskip

In the present paper we give an elementary proof that $R_{H,K}$ given by \eqref{e:bifBm} is nonnegative definite  for $H>1$ and $0<2HK\le 1$, which implies existence of bifractional Brownian motion in this case. (Our proof also works in the case $0<H\le 1$ for $0<K\le 1$ and $0<2HK\le 1$). Moreover, we show that for these parameters $R_{H,K}$  can be written as a sum of two nonnegative definite functions, hence bifBm is equivalent in law to a sum of two independent Gaussian processes. Our theorem is the following:
\begin{thm}
\label{thm:main}
\begin{itemize}
\item[(i)] For any $0<\gamma\le 1$ the function
 \begin{equation}
  C_{\gamma}(s,t)=(s+t)^\gamma-\left(max (s,t)\right)^\gamma, \qquad s,t\ge 0
\label{e:C}
  \end{equation}
is nonnegative definite.
\item[(ii)] Suppose that $\gamma>0$. Then the function
 \begin{equation}
  Q_{\gamma}(s,t)=\left(max (s,t)\right)^\gamma-\abs{t-s}^\gamma, \qquad s,t\ge 0
\label{e:Q}
  \end{equation}
is nonnegative definite if and only if $\gamma\le 1$.
\item[(iii)] Let $H>0$ and $0< K\le 1$ be such that $2HK\le 1$ and let $R_{H,K}$ be defined by \eqref{e:bifBm}, then $R_{H,K}$ is nonnegative definite. Moreover
\begin{equation}
 \label{e:decomp}
 R_{H,K}(s,t)=2^{-K}C_K(s^{2H}, t^{2H})+2^{-K}Q_{2HK}(s,t).
\end{equation}
\end{itemize}
\end{thm}
The proof of this theorem will be given in the next section.
Note that if $H>1$, then the condition $2HK\le 1$ clearly implies that  $0<K \le 1$ authomatically holds, hence \textit{(iii)} gives existence of bifBm in this case.

Theorem \ref{thm:main} \textit{(ii)} leads to the following decomposition of the standard fractional Brownian motion with Hurst parameter $H<\frac  12$ (see \eqref{e:fBm}). It has a very simple nature, but we were not able to find it in literature.

\begin{wn}\label{wn:1.4}
 Suppose that  $0<H\le \frac 12$. Assume that  $\zeta$ is a centered Gaussian process with covariance $Q_{2H}$ given by \eqref{e:Q} and  $\beta$ is a standard Brownian motion independent of $\zeta$  then 
 the process
 \begin{equation*}
  \frac 1{\sqrt 2}\bigg(\zeta(t)+\beta(t^{2H})\bigg), \qquad t\ge 0
 \end{equation*}
is a fractional Brownian motion with covariance \eqref{e:fBm}.
\end{wn}

This follows immediately from  Theorem \ref{thm:main} (ii) since
 \begin{equation*}
  S_H(s,t)=\frac 12\left( Q_{2H}(s,t)+(\min(s,t))^{2H}\right).
 \end{equation*}
  and $2H\le 1$.

\begin{rem} a)
Observe that the fact that  $C_\gamma$ given by \eqref{e:C} is nonnegative definite is not new, since a centered Gaussian process with this covariance appeared in \cite{DurieuWang} as a limit process in an infinite urn scheme. Our proof of positive definiteness of $C_\gamma$ is different, and it is direct.

 b) After having written the present note, we found that the covariance function $Q_\gamma$ was mentioned without a proof on the webpage 
 {\small\texttt{ https://mathoverflow.net} \texttt{/questions/44528/a-simple-decomposition-for-fractional-brownian-motion} \texttt{-with-parameter-h1-2}}, along with the decomposition of fractional Brownian motion given in Corollary \ref{wn:1.4}.
 We were not able to find any subsequent publications.
 More recently, we also learned that this decomposition of fBm was discovered independently by Ivan Nourdin (private communication), but never published. In both cases 
 it was mentioned that the process with covariance  $Q_\gamma$ appeared in the context of Gaussian random fields.
  To the best of our knowledge, our proof of positive definiteness of $Q_\gamma$ for $\gamma \le 1$, is new, as well as the necessary condition.
 
 c) Now that we have shown that $R_{H,K}$ is nonnegative definite for $H>1$ and $2HK\le 1$ we clearly  also have Lei-Nualart decomposition of the covariance of  fractional Brownian motion with Hurst parameter $HK$ in terms of the covariance bifractional Brownian motion and time rescaled Lei-Nualart process (see \cite{LeiNualart}):
 \begin{equation}
  S_{HK}(t,s)=2^{K-1}R_{H,K}(t,s)+\frac 12\left( s^{2HK}+ t^{2HK}-(t^{2H}+s^{2H})^K \right) .
\label{e:LeiNualart}
  \end{equation}
The second term is positive definite for $K\le 1$.
\end{rem}

Unfortunately, Theorem \ref{thm:main} probably does not cover the whole range of parameters for which \eqref{e:bifBm} is a covariance function. The simulations of \cite{LifshitsVolkova} suggest that for any $H>1$ there exists some $\bar K_H$ satisfying $\frac 1{2H}< \bar K_H<\frac 1H$ such that  $R_{H,K}$ is nonnegative definite also for any $0<K<\bar K_H$. This has not been formally proved.

However, for large $H$ the bound $K\le \frac 1{2H}$ is not far short of the optimal one, since we have the following:
\begin{prop}
\label{prop:largeK}
If a sequence $(H_n, K_n)_{n=1}^\infty$ is such that   $R_{H_n,K_n}$ are nonnegative definite and 
\begin{equation*}
 \lim_{n\to \infty} H_n=+\infty,
\end{equation*}
then
\begin{equation}
\limsup_{n\to \infty}2 H_nK_n\le 1.
\end{equation}
\end{prop}
Note that there is some discrepancy between this result and the simulations in \cite{LifshitsVolkova} for large $H$. The latter  suggested that the critical value of $HK$ should be approaching $0.6$ as $H\to \infty$.

\medskip

In the following proposition we collect some simple properties of   
the newly defined bifractional Brownian motion with $H>1$. They are analogous to the ones obtained for $H<1$ and the proofs are essentially the same as for $H<1$.

\begin{prop}\label{prop:properties}
 Suppose that $\xi$ is a centered Gaussian process with covariance function \eqref{e:bifBm} and $H>1$, $2HK\le 1$. 
 \begin{itemize}
  \item[(i)]\begin{equation*}
             2^{-K}\abs{t-s}^{2HK}\le E\abs{\xi(t)-\xi(s)}^2\le 2^{1-K}\abs{t-s}^{2HK}, \qquad s,t\ge 0.
            \end{equation*}
\item[(ii)] $\xi$ has a  modification whose trajectories are locally H\"older continuous with any index $\kappa\in(0,HK)$.
\item[(iii)] The continuous modification given by (ii) is not a semimartingale for $2HK<2$.
\item[(iv)] The process $(2^{(k-1)/2}(\xi_{T+t}-\xi_T))_{t\ge 0}$, converges in the sense of finite dimensional distributions, as $T\to \infty$, to  fractional Brownian motion with Hurst parameter $HK$. If additionally $\xi$ has continuous sample paths, then the convergence is in law in the space of continuous functions on any interval $[0,M]$, $M>0$.
 \end{itemize}
\end{prop}

Property (i) means that bifractional Brownian motion is a quasihelix in the sense of Kahane \cite{Kahane}.

For $H<1$ the bifractional Brownian motion is not a semimartingale also for $HK=1$ (\cite{RussoTudor}). We suspect that this is also in the case $H>1$, but it would require a different proof.

 \section{Proofs}
We use the notation:
 \begin{equation}
  s\wedge t=\min(s,t), \qquad s\vee t=\max(s,t).
 \label{e:min}
 \end{equation}
 \proofof{Theorem \ref{thm:main}}
 \textit{(i)} If $\gamma =1$, then $C_1(s,t)=s\wedge t$, hence it corresponds to Brownian motion.
 
 Now assume that $0<\gamma<1$. We can write
 \begin{equation*}
  C_\gamma(s,t)=\gamma\int_0^{s\wedge t}((s\vee t)+u)^{\gamma-1}du.
 \end{equation*}
By assumption $\gamma-1<0$ and we have
\begin{align*}
  C_\gamma(s,t)=&\gamma\int_0^{s\wedge t}(s+u)^{\gamma-1}\wedge (t+u)^{\gamma-1}du\notag\\
  =&\gamma \int_0^\infty \int_0^\infty \ind{u\le s}\ind{u\le t}\ind{r\le (s+u)^{\gamma-1}}
  \ind{r\le (t+u)^{\gamma-1}}
  &\label{e:1}drdu.
\end{align*}
Here $\indn$ denotes the indicator function.
Thus, from the last equality it follows that for any $n\in \N$, $t_1,\ldots, t_n\ge 0$ and $z_1,\ldots, z_n\in \C$ (complex numbers) we have
\begin{equation*}
 \sum_{j,k=1}^n z_jC_\gamma(t_j,t_k)\overline{z_k}=\gamma \int_0^\infty \int_0^\infty
\abs{ \sum_{j=1}^n z_j\ind{u\le t_j}\ind{r\le (t_j+u)^{\gamma-1}}}^2\hskip -3pt drdu\ge 0.
\end{equation*}
This finishes the proof of part $(i)$.

\medskip
Part $(ii)$:
If $\gamma=1$, then $Q_1(s,t)=s\wedge t$ is the covariance function of Brownian motion.

Now suppose that $0<\gamma<1$. We use the fact that for $0<\gamma<1$ and  $x\ge 0$ we have  
\begin{equation}
 x^\gamma =c_\gamma \int_0^\infty\frac  {1-e^{-xy}}{y^{\gamma+1}}dy,
 \label{e:x_gamma}
\end{equation}
where 
\begin{equation*}
 c_\gamma=\left(\int_0^\infty \frac {1-e^{-y}}{y^{\gamma+1}}dy\right)^{-1}=\frac \gamma{\Gamma(1-\gamma)}>0.
\end{equation*}
Here $\Gamma$ denotes the usual Euler Gamma function.

We write $Q_\gamma$ with help of \eqref{e:x_gamma} and use $t\vee s-\abs{t-s}=s\wedge t$ obtaining
\begin{align}
 Q_\gamma(s,t)=&c_\gamma \int_0^\infty \frac {e^{-\abs{t-s}y}-e^{-(s\vee t)y}}{y^{1+\gamma}}dy\notag\\
 =&c_\gamma \int_0^\infty \frac {e^{-\abs{t-s}y}(1-e^{-(s\wedge t)y})}{y^{1+\gamma}}dy.\label{e:2.3}
 \end{align}
 Note that for $y>0$
 \begin{equation*}
 e^{-\abs{t-s}y}=Ee^{iy(\eta_t-\eta_s)}
 \end{equation*}
where $\eta $ is the standard Cauchy process, that is, a process with stationary independent increments and such that $\eta_1$ has the standard Cauchy distribution.
 We also use 
 \begin{equation*}
  \frac{1-e^{-(s\wedge t)y}}y=\int_0^\infty e^{-ry}\ind{r\le s}\ind{r\le t}dr.
 \end{equation*}
(Note that the above can be also represented as the covariance function of the process $\int_0^t e^{-ry/2}dW_r$, where $W$ is a Brownian motion.)
Hence from \eqref{e:2.3} we obtain
 \begin{equation*}
  Q_\gamma(s,t)
 = c_\gamma \int_0^\infty \int_0^\infty  y^{-\gamma}e^{-y r} 
 E e^{iy(\eta_t-\eta_s)} \ind{ r\le t}\ind{r\le s}dr dy.
\end{equation*}
Consequently,  for any $n\in \N$, $t_1,\ldots, t_n\ge 0$ and $z_1,\ldots, z_n\in \C$ we have
\begin{equation*}
 \sum_{j,k=1}^n z_jQ_\gamma(t_j,t_k)\overline{z_k}=c_\gamma  \int_0^\infty \int_0^\infty y^{-\gamma} e^{-yr}
E\abs{ \sum_{j=1}^n z_j\ind{r\le t_j} e^{iy\eta_{t_j}}}^2 drdu\ge 0.
\end{equation*}
Hence $Q_\gamma$ is nonnegative definite for $0<\gamma<1$.

\medskip

Now suppose that $\gamma>1$.
If $Q_\gamma$ were nonnegative definite, then the function 
\begin{equation*}
 f(a):=Q_{\gamma}(1,1)+Q_{\gamma}(1+a,1+a)-2Q_\gamma(1,1+a)
 =1+2a^\gamma -(1+a)^\gamma, \qquad a\ge 0
\end{equation*}
would have to be nonnegative for all $a\ge 0$. Note however, that $f(0)=0$ and
\begin{equation*}
 f'(a)=2\gamma a^{\gamma-1}-\gamma (1+a)^{\gamma-1}.
\end{equation*}
Hence, if $\gamma >1$ the function $f'$ is negative in some interval $(0,\varepsilon)$ with $\varepsilon>0$ and therefore  on this interval $f$ is decreasing and takes negative values.
This finishes the proof of \textit{(ii)}.

\textit{(iii)}: It is clear that  \eqref{e:decomp} holds. Hence \textit{(iii)} follows direcly from \textit{(i)} and \textit{(ii)}, since we have assumed that  $0<K<1$ and $0<2HK<1$.
\qed

\proofof{Proposition \ref{prop:largeK}} Without loss of generality we can consider a sequence $K_n,H_n$ satisfying the assumptions of the proposition and such that 
$\lim_n{2H_nK_n}$ exists
\begin{equation*}
 \lim_{n\to \infty}2H_nK_n=\gamma.
\end{equation*}
Clearly $K_n\to 0$ and  $\gamma\le 2$ since $K\le \frac 1H$ is a necessary condition for $R_{H,K}$ to be nonnegative definite. We have to show that $\gamma\le 1$.

Suppose that $0\le s\le t$, and $t>0$. Then we can write
\begin{equation*}
 R_{H,K}(s,t)=2^{-K}\left[t^{2HK} \left(1+ (\frac{s}{t})^{2H}\right)^K-\abs{t-s}^{2HK}\right], 
\end{equation*}
from which it follows immediately that under our assumptions 
\begin{equation*}
\lim_{n\to\infty} R_{H_n,K_n}(s,t)=t^\gamma-\abs{t-s}^\gamma
\end{equation*}
Hence, in general, for any $s,t\ge 0$
\begin{equation*}
\lim_{n\to\infty} R_{H_n,K_n}(s,t)=(s\vee t)^\gamma-\abs{t-s}^\gamma=Q_\gamma(s,t).
\end{equation*}
By the ``only if'' part of Theorem \ref{thm:main} (ii),
if $\gamma>1$ then $Q_\gamma$ is not nonnegative definite, 
which contradicts the assumption that $R_{H_n,K_n}$ were nonnegative definite.
\qed

\proofof{Proposition \ref{prop:properties}}
(i) The proof is the same as that of Proposition 3.1 in \cite{HoudreVilla}.\\
(ii) Follows direcly from (i) and the Kolmogorov criterion, and the fact that the process is Gaussian.\\
(iii) Follows from (i) and Lemma 2.1 in \cite{BGT_fBm_density}. In fact, from the proof of that lemma we have that $\xi$ has infinite quadratic variation.\\
(iv) To obtain the convergence of finite dimensional distributions it suffices to adapt the proof of Theorem 2.1 in \cite{MaejimaTudor}, since it only used that $K<1$ and $HK<1$. Tightness in the continuous case follows from part (i) and Gaussianity of the process.

\medskip
\textbf{Acknowledgement.} The author thanks M. Lifshits for pointing to reference \cite{DurieuWang} and T. Bojdecki and  L.G. Gorostiza for their  comments.

 \end{document}